\documentclass[12pt,leqno]{article}
\usepackage{amssymb}
\newcommand{\dd}{{\rm \kern 3pt I\kern-9pt d}}

\newcommand{\Abar}{{\backslash\kern-8pt A}}

\font\helv Helvetica at 14pt
\newcommand{\CC}{\mbox{{\helv \i}{\rm\kern-5pt C}}}

\begin{document}
\begin{center}
\LARGE
Stochastic approach for the subordination in Bochner sense
\end{center}

\hspace{10cm} Nicolas Bouleau

\hspace{10cm} Ecole des Ponts, Paris\\

It is possible to construct (cf [1] and [2]) a positive process $(Y_{\psi,t})$ indexed by Bernstein functions and time which, for fixed $t$, is Markovian with respect to the composition of Bernstein functions, and for fixed $\psi$, is the subordinator associated with the Bernstein function $\psi$. In the same manner a realization $(X_{Y_{\psi,t}})$ can be obtained of all the subordinated processes of a Markov process $(X_t)$.

This probabilistic interpretation of the initial idea of Bochner allows to construct martingales and to apply stochastic calculus to questions related to symbolic calculus for operators semi-groups.

We study here a particular branch of the subordination process : The homographic branch which has the advantage of being connected with some works on positive infinitely divisible diffusions and brings on them a different point of view. Nevertheles, what follows would apply, except technical difficulties, to any branch of the subordination process.

It is easy to see that the functions
$$f_a(x)=\frac{xe^a}{1+x(e^a-1)},\quad a\geq 0,x\geq 0$$
satisfy $f_a\circ f_b=f_{a+b}$ and are Bernstein functions (cf [3]) associated with the subordinators 
\begin{equation}
Y_t^a=(e^a-1)\sum_{k=1}^{N_{t/(1-e^a)}}\;E_k
\end{equation}
where $(N_t)$ is a standard Poisson process and the $E_k$'s are i.i.d. exponential random variables independent of $N$. The L\'evy measure of $(Y_{a,t})$ is therefore
$$\nu_a(dy)=\frac{1}{(e^a-1)^2}e^{-y/(e^a-1)}\, dy\quad\mbox{ on } \mathbb{R}_+.$$
The relation $P_au(x)=\mathbb{E}[u(Y_x^a)]$ for $u : \mathbb{R}_+\mapsto\mathbb{R}_+$ defines a Markovian semi-group $(P_a)_{a\geq 0}$ with generator $Au(x)= xu^\prime(x)+xu^{\prime\prime}(x)$ which is the transition semi-group of the diffusion $(Z_a)_{a\geq 0}$ solution to the sde
\begin{equation}
dZ_a=\sqrt{2Z_a}dB_a+Z_ada,\quad Z_0=z\geq 0.
\end{equation}
We now define a two parameters process $(Y_{a,t})_{a\geq 0, t\geq 0}$ by chosing for $a_1<a_2\cdots<a_n$ the joint law of the processes
$$(Y_{a_1,t})_{t\geq 0},\ldots, (Y_{a_n,t})_{t\geq 0}$$
to be that of the process
$$\left(Y_t^{a_1},Y_{Y_t^{a_1}}^{a_2-a_1},\ldots,Y^{a_n-a_{n-1}}_{Y^{a_{n-1}-a_{n-2}}_{\ddots Y_t^{a_1}}}\right)_{t\geq 0}$$
where $(Y_t^{a_1})_{t\geq 0},(Y_t^{a_2-a_{1}})_{t\geq 0},\ldots,(Y_t^{a_n-a_{n-1}})_{t\geq 0}$ are independent subordinators of type (1).

A version of the process $(Y_{a,t})$ may be chosen such that $Y_t=Y_{.,t}$ be right continuous and increasing with values in $\mathcal{C}(\mathbb{R}_+, \mathbb{R}_+)$ and with independent stationary increments, and for fixed $t$ $(Y_{a,t})_{a\geq 0}$ has the same law as $(Z_a)$ for $Z_0=t$.

Because of the formula $\nu_aP_b=\nu_{a+b}$, we can construct the process $(Z_a)_{a>0}$ with the entrance rule $(\nu_a)$ what defines a positive $\sigma$-finite measure $m$ on the space $\mathcal{C}(\mathbb{R}_+^\ast, \mathbb{R}_+)$ as the `law' of the process $(Z_a)$ under the entrance law $(\nu_a)$.

This measure $m$ is the L\'evy measure of the process $(Y_t)_{t\geq 0}$. This can be seen in the following way. Let $\mu$ be a positive measure with compact support on $\mathbb{R}_+$, let us put $<Y_t,\mu>=\int Y_{a,t}\,\mu(da)$. The relation
\begin{equation}
\mathbb{E}e^{-<Y_t,\mu>}=\exp{\int(e^{\int Z_\alpha\mu(d\alpha)}-1)\,dm}
\end{equation}
is easy to prove when $\mu$ is a weighted sum of Dirac masses by the computation
$$\mathbb{E}e^{-<Y_t,\lambda\varepsilon_a>}=\mathbb{E}e^{-\lambda Y_{a,t}}=e^{-tf_a(\lambda)}=e^{t\int(e^{-\lambda y}-1)\;\nu_a(dy)}$$
noting that $\nu_a$ is the law of $Z_a$ under $m$. And for general $\mu$ (3) is obtained by weak limit.

The law of the subordinators $<Y_t,\mu>$, may be studied by using the fact that if we decompose $\mu$ into
$$\mu=1_{[0,x]}.\mu+1_{(x,\infty)}.\mu=\mu_1+\mu_2$$ the following representation holds
$$<Y_t,\mu>=<Y_t,\mu_1>+\int Y^{b-x}_{Y_{x,t}}\,\mu_2(db)$$
where $(Y_t^\beta)_{\beta\geq 0,t\geq 0}$ is independent of $(Y_{a,t}, a\leq x, t\geq 0)$.

One obtains for $\lambda\geq 0$, 
$$\mathbb{E}\exp\{-\lambda<Y_t,\mu>\} =\exp\{-t(\frac{1}{2}-g_x^\prime(0,\lambda))\}$$
where $g$ is the positive solution decreasing in $x$ of 
$$\left\{
\begin{array}{rl}
g_{x^2}^{\prime\prime}=&g.(\frac{1}{4}+\lambda.\mu)\\
g(0,\lambda)=&1
\end{array}
\right.
$$ and if  we put $\xi(x,\lambda)=\frac{1}{2}-\frac{g_x^\prime(x,\lambda)}{g(x,\lambda)}$, the process
$$M_{a,t}=\exp\{-\xi(a,\lambda)Y_{a,t}+\xi(a,\lambda)t-\lambda\int_0^aY_{\alpha,t}\,\mu(d\alpha)\}$$
is a two parameters martingale for the filtration of $(Y_{a,t})$.

A similar study may be performed with the family $g_a(x)=\frac{x}{1+ax}$ which corresponds to the diffusion 
\begin{equation}
Z_a=z+\int_0^a\sqrt{2Z_b}\,dB_b.
\end{equation} wich is the square of a Bessel process of exponent 0, cf [4].

These results were obtained by Pitman and Yor [5], [6], by interpretating the process $Y_{a,t}$ of the case (4) as $\ell_{\tau_t}^a$ where $\ell_t^a$ is the family of brownian local times and $\tau_t=\inf\{s:\ell^0_s=t\}$, and using the Ray-Knight theorems. In this framework the measure $m$ appears to be the image of the Ito measure of the Brownian excursions.\\

\noindent References

[1] {\sc N. Bouleau} and {\sc O. Chateau}  ``Le processus de la subordination", {\it C. R. Acad. Sc. Paris} t. 309, s1, 625-628, (1989).

[2] {\sc O. Chateau} {\it Quelques remarques sur les processus \`a accroissements ind\'ependants et la subordination au sens de Bochner} Th\`ese Univ. Paris VI, (1990).

[3] {\sc C. Berg} and {\sc G. Forst}, {\it Potential theory on locally compact groups}, Springer (1975).

[4] {\sc T. Shiga} and {\sc S. Watanabe}, ``Bessel diffusions as a one parameter family of diffusion processes" {\it Z. f. Wahr.} 27, 37-46, (1973).

[5] {\sc J. Pitman} and {\sc M. Yor}, ``A decomposition of Bessel bridges",{\it Z. f. Wahr.} 59, 425-457, (1982).

[6] {\sc J. Pitman} and {\sc M. Yor}, ``Sur une d\'ecomposition des ponts de Bessel", Lect. N. in Math. 923, p.276 {\it et seq.}, (1982).
\end{document}